\documentclass[leqno, draft]{amsart}
\usepackage[c5paper, text={126truemm,178truemm},centering]{geometry}
\usepackage{xcolor}
\usepackage{tikz}
\usepackage{comment}
\usepackage{pict2e}
\usepackage[utf8]{inputenc}
\usepackage{multicol}
\usetikzlibrary{calc,arrows,decorations.markings,intersections,matrix}
\newcommand\includetikz[2][x=1mm,y=1mm]{%
 \IfFileExists{#2.tikz}{%
	\begin{tikzpicture}[#1]\node at (0,0) {\input{#2.tikz}};\end{tikzpicture}%
 }{%
	\@latex@error{No usable file `#2.tikz' can be found}%
	 {I could not locate the file...}%
 }
}
\def\Bf#1{\ifmmode\boldsymbol{#1}\else{\rmfamily\bfseries#1}\fi}
\def\<{\langle}
\def\>{\rangle}

\usepackage[backref=page]{hyperref}
\hypersetup{colorlinks, citecolor=[rgb]{0.0,0.3,0.1}, linkcolor=[rgb]{0.3,0.0,0.1},
						urlcolor=[rgb]{0.0,0.0,0.4}}
\makeatletter
\def\BR@backref#1{{\upshape
 \begingroup
  \csname @safe@activestrue\endcsname
  \expandafter\providecommand\csname brc@#1\endcsname{0}%
  \expandafter\providecommand\csname brcd@#1\endcsname{0}%
  \csname @safe@activesfalse\expandafter\endcsname
  \ifBR@BackrefAlt
   \ifx\backrefentrycount\BR@BackrefEntryCountUnused
   \else\BR@PopulateEntryCount{#1}\fi
   \expandafter\backrefalt\csname brc@#1\expandafter\endcsname
   \csname brl@#1\expandafter\endcsname
   \csname brcd@#1\expandafter\endcsname
   \csname brld@#1\endcsname
  \else
   \expandafter\backref\csname br@#1\expandafter\endcsname
  \fi
 \endgroup
 \par
}}
\makeatother
\renewcommand*{\backref}[1]{\ifx#1\relax\relax\else{\rm~$\langle$#1$\rangle$}\fi}

\theoremstyle{plain}
\newtheorem{theorem}{Theorem}[section]

\newtheorem{question}[theorem]{Question}
\newtheorem{statement}[theorem]{Statement}

\theoremstyle{definition}

\newtheorem{remark}[theorem]{Remark}

\numberwithin{figure}{section}

\title[A local support theorem]{                   
					Local support theorem for the exponential Radon transform
}

\author[E. Dinnyés, T. Ódor]{Enikő Dinnyés${}^{*,\#}$ \and Tibor Ódor${}^{**,\#}$}
\thanks{${}^*$ Department of Probability and Statistics, Eötvös Loránd University, Budapest}
\thanks{${}^{**}$ Department of Geometry, Bolyai Institute, University of Szeged}
\thanks{${}^{\#}$ Rényi Institute of Mathematics, Budapest}
\thanks{This work on the project leading to this application has received
funding from the European Research Council (ERC) under the European Union’s Horizon 2020 research and innovation programme (grant agreement No. 741420).}

\begin{document}


\begin{abstract}
We prove a local support theorem for the exponential Radon transform for functions of exponential decay at infinity. We also show that our decay condition is essentially sharp for the classical Radon transform for hyperbolic type domains as holes, by showing streched exponential counterexamples. This shows a difference of the support theorems for compact domains, where the decay has to be just faster than any polynomial. Also, gives a refinement for non-compact domains, where support theorem was proved only for functions with compact support. Our method is a version of~\cite{Strichartz} and~\cite{BomanUniqueness}.
\end{abstract}

\maketitle

\section{\label{sec:intro}
                              Introduction
}

\bigskip

\smallskip
Exterior problems are a special type of region-of-interest tomography, where the internal part of the projections are unavailable. This is a version of truncated projections: it falls into the category of limited-angle problems. Unique reconstruction is not always possible (as we approach the internal domain that must be avoided, the angle between possible scans is more and more limited).

Our support theorems speak about this set of problems. We measure some version of the Radon transform outside an internal domain and if we find it to be zero everywhere then we want to conclude that the original function is zero outside the given domain.

In this article we discuss the case when we only have the data yielded by a version of the exponential Radon transform along lines not intersecting a convex but not necessarily compact (not necessarily bounded) closed set $K$. Furthermore, we do not require the density function of interest to have compact support, but of exponential decay.

\smallskip
If the hole $K$ is convex and compact, and we know the integrals over lines avoiding $K$ the support theorem is known for $f$ decaying faster than any polynomial (for each integer $k>0, \;\;|x|^k\,f(x)$ is bounded)~\cite{Helgason}. On functions this is a weaker restriction than exponential 
decay in Theorem~\ref{th:supp}, but $K$ has to be compact, unlike in our theorem. 
On the other hand, if $K$ is non-compact, but closed a similar support theorem was proved for the classical Radon transform for functions of squared exponential decay e.g. in~\cite{Odor}. For functions with compact support there is a general support theorem for generalized analytic Radon 
transforms~\cite{BomanQuinto}.

For hyperbolic type domains (containing two half-lines with the same endpoints), our results in Theorem~\ref{th:supp} are essentially sharp. There is a 
counter-example with decay $\exp(-\mu \cdot x^{1-\epsilon})$ for  any $\mu > 0$ and $0< \epsilon < 1$. 
For parabolic domains (all directed half-lines in $K$ are parallel) we have only counter-example with decay $\exp(-\lambda \cdot x^{1/2-\epsilon})$. It is not known, whether this counter-exeample is essentially sharp, or there is one with faster decay. Note that both hyperbolic and parabolic domains are assumed to be part of an angle, the intersection of two half-planes with different normal vectors.

These results show that there is a trade-off between decay conditions and geometric properties of the hole $K$.

\bigskip
Let $K\subset\mathbb R^2$ be a 
closed convex set. We say that the function $f\colon K^c\to \mathbb R$ has an \emph{exponential decay} if there exists a $\mu > 0$ such that for every $\varepsilon > 0$ the integral $\int_{K^c_\varepsilon} |f(x)|\, e^{\mu |x|}\, dx$ is finite, where $K^c_\varepsilon$ is the set of those points in $\mathbb R^2$ whose distance from $K$ is not less than $\varepsilon$.

\smallskip

\smallskip
The \emph{exponential Radon transform} with parameter $\lambda$ is defined as
$$R_\lambda f(\omega,p) = \int_{\< x,\omega \>=p} e^{-\lambda\< x,\omega^\perp\>} f(x)dx = \int_{-\infty}^\infty e^{-\lambda u} f(p\omega+u\omega^\perp)\,du,$$
along the line $l(\omega, p) := \left\{ x \in\mathbb{R}^2 : \< x,\omega\> = p \right\}$. Here $dx$ is the Lebesgue measure on the line $\< x, \omega\> = p$, and $\omega^\perp$ is 
$(\cos \omega, \sin\omega)^\perp = (\cos (\omega + \frac{\pi}{2}), \sin(\omega + \frac{\pi}{2})) = (-\sin\omega, \cos\omega)$. 
With $\lambda=0$ we get the classical Radon transform of $f$.

Note that the line $l(\omega, p)$ is the same set as $l(-\omega, -p)$, but admits different orientation. Also note that although defined on the same set, the integral $(R_\lambda f)(\omega, p)$ usually differs from $(R_\lambda f)(-\omega, -p)$, unlike the special case of the classical Radon transform where $Rf(\omega,p)=Rf(-\omega,-p)$.

\bigskip
Our goal is to prove the following support theorem for a less restrictive class of functions than the functions of exponential decay, for the exponential Radon transform. But we state it in this simple form first, with its proof in the next section.

\begin{theorem}\label{th:supp}
Let $K\subset\mathbb R^2$ be a (not necessarily bounded) closed convex set not containing a whole line, and $f\colon K^c\to \mathbb R$ be a locally integrable function of exponential  decay. Assume that the exponential Radon transform $R_\lambda f(\omega, p) = 0$ along every line $l(\omega, p)$
not intersecting $K$. Then $f = 0$ in $L^1_{loc}(K^c)$.
\end{theorem}

We will prove this theorem as a consequence of Theorem~\ref{th:supp2}.

\section{\label{sec:motiv}
                              Motivation: 
                              counterexamples for streched exponential decay
}

The stretched exponential function $e^{-z^{\beta}}$ for a $\lambda > 0$ 
and $0 < \beta < 1$ is obtained by inserting a fractional power into the exponential 
function.

The following two counter-examples work for classical Radon-transform, that is for $\lambda = 0$. 

For a parabolic convex closed set $K$, that every oriented half-line is parallel, a simple example is 
$$f(z)=e^{-z^{1/2-\epsilon}}, \quad z\in\mathbb C\;,$$
where $0< \epsilon < 1/2$ is arbitrary.
This is a complex analytic function apart from the half-line $$\mathbb R_-=\{z: \mathrm{Im}z=0, \mathrm{Re}z\le 0\}\;;$$ it is important that the real part of $z^{1/2-\epsilon}$ is positive because its angle is a bit less than half the original angle, and $|f(z)|=e^{-Re(z^{1/2-\epsilon})}$ tends to zero as 
$e^{-\theta}_\epsilon \cdot r^{1/2-\epsilon}$ for some $\theta_\epsilon > 0$ when $r = |z|\rightarrow\infty$ so that $z$ moves along a straight line, e.g. if we put $z=az_0$ with $z_0\in\mathbb{C}, a\in\mathbb{R}, a>0$ and $a\rightarrow\infty$.

$f(z)$ being complex analytic means its complex integral along any closed curve not intersecting $\mathbb R_-$ is zero. Fix a straight line not intersecting $\mathbb R_-$, fix a point on the line, and draw a half circle from this point on that side of the line that does not contain $\mathbb R_-$. Let the closed curve be this half circle with its diameter along the straight line. The complex integral of $e^{-z^{1/2-\epsilon}}$ is zero along this curve. Take the limit of this integral as the radius tends to infinity. The complex integral along the straight line is a constant multiple of the combination of the real integrals consisting of the real part and the imaginary part separately; these tend to the Radon transform of the real part and the Radon transform of the imaginary part respectively; and the norm of the integral along the half circle tends to zero. So the Radon transforms of the real functions $Re(f)$ and $Im(f)$ along any straight line not intersecting $\mathbb R_-$ must be zero. 

In fact, this example is applicable to any domain $K$ containing a half line.

\medskip
The following example for hyperbolic domains works in a similar way. 

Let 
$$f(z)=\exp(-z^{1-\varepsilon}), \quad z\in\mathbb C$$
for some $0<\varepsilon<1$ of our choice. The same technique is applicable where the real part of $z^{1-\varepsilon}$ is positive: in the concave angular domain $\left(-\frac{1}{(1-\varepsilon)}\frac{\pi}{2},\;\frac{1}{(1-\varepsilon)}\frac{\pi}{2}\right),$ and $K$ must contain all the rest. 

When $\varepsilon$ is close to 0, the domain to be covered by $K$ is a wide angular domain, close to the $\mathrm{Re}z > 0$ half plane. But applying an affine transformation for a general $K$ containing any angular domain, mapping it to the required one, this second example can be applied in the case of any non-bounded closed convex hyperbolic $K$.

\medskip
Assume that $K\subset {\mathbb R^2}$ is closed, convex, $int(K)\neq \emptyset$, unbounded  and does not contain a full line. 
Then for every $x\in \partial K$ we have a half line $h\subset K$ with endpoint $x$, and its relative interior in $int(K)$.
(Let $y\in \partial{K}$, and $y\rightarrow\infty$ along the border of $K$. Then the limit of $xy$ is $h$. If $h\subset \partial K$ then move $y$ to the other direction along the border of $K$. In the limit, we get the half-line $g$. If $g\subset \partial K$ then $K$ is an angular domain, for which the statement is obviously true. Otherwise $h$ or $g$ is a suitable half-line. Assume it is $h$.)
Let $s$ be a support line of $K$ in $x$. Then apply an affine transform $A$ which maps $h$ to the positive half of the $x_1$ axis, and $s$ to the $x_2$ axis. 

Let $f$ be as before, but restricted to $AK^c$. Then $R_\lambda f(\omega,p)=0$ for every line $l(\omega, p)$ that avoids $AK$. 

Let $f_A(x)=f(Ax)$. If $l(\omega, p)$ avoids $K$ then $Al(\omega,p)=:l(\omega^\star, p^\star)$ avoids $AK$, and vice versa. 

If $R_\lambda f(\omega^\star,p^\star)=0$ then $R_\lambda f_A(\omega,p)=0$, and vice versa.

This yields a function whose exponential Radon transform $R_\lambda$ is zero on lines avoiding $K$ and which decays streched exponentially at the infinity.

By applying convolution, we can make it smooth, even analytic. We leave the details to the reader, which is a standard computation, using that the transform $R_\lambda$ commutes with convolution.

\medskip
These counter-examples show that if $K$ is non-bounded then its behaviour is rather different from the compact case, where we have support theorems 
for functions decaying faster than any polynomial~\cite{Helgason},~\cite{Strichartz}. 

For functions with compact support similar theorems were known 
for a general $K$~\cite{BomanQuinto}. 
The novelty in this paper is extending the appropriate function class, and showing that less than exponential decay is not enough.

Theorem~\ref{th:supp} is a support theorem for non-bounded $K$.

For a general exponent $\lambda > 0$, the functions have to decay faster than $e^{-(\lambda t + \epsilon |t|)}$ for any fixed $0 < \epsilon < \lambda$. Note, that in the exponent the $\lambda t$ term can be negative for negative $t$. Theorem~\ref{th:supp} and Theorem~\ref{th:supp2} shows that exponential decay is enough even for non-bounded $K$ if it does not contain a whole line (that is a trivial condition because if $K$ contained a whole line, only lines parallel to it could avoid $K$).

In case of the classical Radon transform (that is $\lambda = 0$) for parabolic domains still there is a gap in our knowledge. This gap cannot be directly attacked by our technique. 

\begin{question} 
What are the close to exact exponential type decay conditions for parabolic convex closed sets $K$? 
\end{question}

\section{\label{sec:exp}
                              A local support theorem for the exponential Radon transform for functions of exponential decay 
}

Consider the following conditions as local conditions.

Let $(\omega_0, p_0)\in \mathbb S^1\times \mathbb R$ be fixed,  $p_0>0$. Let 
$\Omega \subset \mathbb S^1$ be an open subset of  
$\mathbb S^1$ containing $\omega_0$. Let $\mathcal L(\Omega, p_0)$ be the union of all lines $l(\omega, p)$, as point sets, such that $p > p_0$ and $\omega \in \Omega$. 

The complement $K^c$ of any closed convex set $K$ can be constructed as the union of sets belonging to this type with some $\Omega$ and $p_0$, if there is no complete straight line in the border of $K$ (that means $K$ is neither a half plane nor a strip). 
 E.g. the complement of a parabola can be constructed as the union of infinitely many (countably many) of these sets, so in that case we need the below statement for infinitely many $\omega_0$'s and $p_0$'s.
But the conditions for the function $f$ can be much more easily written down for each $\mathcal L(\Omega, p_0)$ separately, that is why we use this description.

We define the space of functions $\mathcal E(\omega_0, \Omega,  p_0)$ as follows. $f\in\mathcal E(\omega_0, \Omega, p_0)$ if and only if

a) $f$ is continuously differentiable on $\mathcal L(\Omega, p_0)$; 

b) $f(x) =f(p\, \omega + u\, \omega^\perp)$ decays locally uniformly faster than every polynomial on every line $l(\omega, p)$ in $\mathcal L(\Omega, p_0)$ $-$ uniformly in $\omega$ $-$,
that means for every $p$ and every $k\in \mathbb N$ there exists a neighbourhood $(p-\varepsilon_p, p+\varepsilon_p)$ such that
for any $p'\in(p-\varepsilon_p, p+\varepsilon_p)$,
independently of $\omega$, there exists a $c_{pk}$ constant such that 
$$\int_{-\infty}^\infty |u^k \, f(p'\, \omega + u\, \omega^\perp)|\, du < c_{pk} < \infty;$$

c) $\partial_{\omega} f(x) =(\partial_{\omega} f)(p\, \omega + u\, \omega^\perp)$, where $\partial_{\omega}$ means directional derivative in the direction $\omega$, decays locally uniformly faster (uniformly in $\omega$) than every polynomial on every line $l(\omega, p)$ in $\mathcal L(\Omega, p_0)$ 
that means for every $p$ and every $k\in \mathbb N$ there exists a neighbourhood $(p-\varepsilon_p, p+\varepsilon_p)$ such that 
for any $p'\in(p-\varepsilon_p, p+\varepsilon_p)$,
independently of $\omega$, there exists a $C_{pk}$ constant such that 
$$\int_{-\infty}^\infty |u^k \, (\partial_\omega f)(p'\, \omega + u\, \omega^\perp)|\, du < C_{pk} < \infty.$$

d) for every $p> p_0$ we demand that $f(p\, \omega_0 + u\, \omega^\perp_0)$ 
decays exponentially, that is, for some $\mu_p > 0$ we have 
$$\int_{-\infty}^\infty\,|f(p\, \omega_0 + u\, \omega_0^\perp)|\, e^{\mu_p |u|}\,du < K_p < \infty$$
(note that in this condition we demand exponential decay only for directions perpendicular to $\omega_0$, not for all $\omega\in\Omega$); 

e) we assume that there is a sequence $\left\{q_k\right\}_{k\in\mathbb N}, \; p_0 < q_k \leq +\infty$, such that $$R^{(k)}_\lambda f(\omega_0, q_k):=\int_{-\infty}^\infty u^k e^{-\lambda u}\, f(q_k\omega_0 + u\omega_0^\perp) \,du = 0,$$ or in a weaker form, the integral exists, finite and is known, for every $k\geq 0$. 

It is possible that $q_k$ is a finite constant ($=q_0$), implying that we know the function $f$ on the line $l(\omega_0, q_0)$. Or all $q_k$ can be $+\infty$. For a fixed $k$, $q_k=+\infty$ means that there is a sequence $q_{kj}\to +\infty$ such that
$$\limsup_{j\rightarrow\infty}\left(\int_{-\infty}^\infty |u^k e^{-\lambda u}\, f(q_{kj}\, \omega_0 + u\, \omega_0^\perp)|\, du\right)< \infty$$ and $\lim_{j\rightarrow\infty}\left(\int_{-\infty}^\infty u^k e^{-\lambda u}\, f(q_{kj}\, \omega_0 + u\, \omega_0^\perp)\, du\right)$ exists and is known.

No monotonicity or other restriction is posed to the sequence $q_k$ or to the function $f$.

\bigskip
These conditions form a very wide 
class compared to the 'usual' class of functions with compact support, or even the functions of exponential decay. This shows the flexibility of our technique.




\emph{Example:} Condition $e)$ is fulfilled by the following function: $$f(x)=\exp(\<\omega_0,x\>^2)\,\exp(-\<\omega_0^\perp,x\>^2)\,\sin(\<\omega_0,x\>^2).$$

For $\<\omega_0,x\>^2=k\pi$ the term $\sin(\<\omega_0,x\>^2)=0$, so $q_{k}=\sqrt{k\pi}$ will be suitable for all $k$ to fulfill condition $e)$.

\medskip
\begin{theorem}\label{th:supp2} Knowing the exponential Radon transform of a function $f$ belonging to the function class $\mathcal E(\omega_0, \Omega, p_0)$ defined on $\mathcal{L}(\Omega,p_0)$, the function $f$ is determined on $\mathcal{L}(\Omega,p_0)$.
\end{theorem}

\medskip
\begin{remark} The system of conditions in Theorem ~\ref{th:supp2} are mixed, compared to polynomial decay: stronger for some directions, but weaker for other directions (no restriction for decay properties at all).
\end{remark}

\medskip
\emph{Proof of Theorem~\ref{th:supp2}.}
First, we write down an equality for the partial derivative of the exponential Radon transform $R_\lambda f(\omega, p)$ with respect to $p$. Knowing the exponential Radon transform along lines avoiding $K$ this derivative can be calculated.

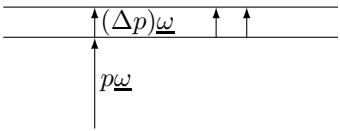
\begin{multicols}{2}

\setlength{\unitlength}{0.4 cm}
\begin{picture}(12,6)(-3.5,-3.5)
\put(-3,0){\line(1,0){11}}
\put(-3,1){\line(1,0){11}}
\put(0,0){\vector(0,1){1}}
\put(0.2, 0.5){\makebox(0,0)[l]{$(\Delta p){\underline\omega}$}}
\put(4,0){\vector(0,1){1}}
\put(5,0){\vector(0,1){1}}
\put(0,-3){\vector(0,1){3}}
\put(0.2,-1.5){\makebox(0,0)[l]{$p{\underline\omega}$}}
\end{picture}

Study the first figure for the partial derivative with respect to $p$.

$\<x,\omega^\perp\>$ does not change along the arrows. The change of $f$ along the short arrows is $\Delta p\,\partial_{\omega} f(x)$, that we divide by $\Delta p$ when differentiating. 

\end{multicols}

So the derivative with respect to $p$ is
$$\partial_p R_\lambda f (\omega, p) = \int_{\<x,\omega\>=p} e^{-\lambda\< x,\omega^\perp\>}\partial_\omega f(x)\, dx,$$
where $\partial_\omega$ means directional derivative in the direction $\omega$.

\medskip
Now let's find an equality for the partial derivative of $R_\lambda f(\omega, p)$ with respect to the angle $\omega$.

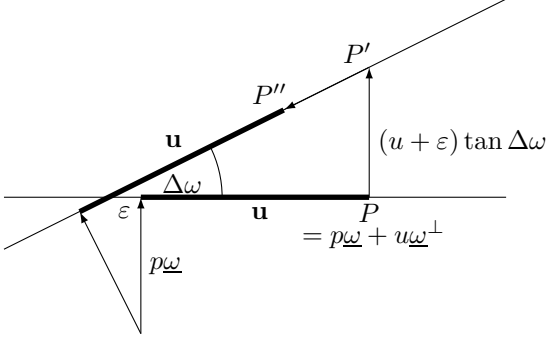
\setlength{\unitlength}{0.6 cm}
\begin{picture}(12,8)(-3.5,-3.5)
\put(-3,0){\line(1,0){11}}
\put(-3,-1.14){\line(2,1){11}}
\put(5,0){\vector(0,1){2.86}}
\put(5, 2.86){\vector(-2,-1){1.87}}
\put(-0.72,0){\arc[0,27]{2.5}}
\put(0.9,0.3){\makebox(0,0){$\Delta\omega$}}
\put(2.6,-0.35){\makebox(0,0){\textbf{u}}}
\put(5.2,1.2){\makebox(0,0)[l]{$(u+\varepsilon)\tan\Delta\omega$}}
\put(-0.36,-0.2){\makebox(0,0)[t]{$\varepsilon$}}
\put(5, -0.15){\makebox(0,0)[t]{$P$}}
\put(5.1, -0.5){\makebox(0,0)[t]{$= p{\underline\omega}+u{\underline\omega}^\perp$}}
\put(5.05, 2.95){\makebox(0,0)[br]{$P'$}}
\put(3.2, 2.05){\makebox(0,0)[br]{$P''$}}
\put(0.2,-1.5){\makebox(0,0)[l]{$p{\underline\omega}$}}
\put(0,-3){\vector(0,1){3}}
\put(0,-3){\vector(-1,2){1.34}}
\linethickness{2pt}
\put(0,0){\line(1,0){5}}
\put(-1.34,-0.32){\line(2,1){4.47}}
\put(0.9, 1.1){\makebox(0,0)[br]{\textbf{u}}}
\end{picture}

We consider the change of the integrand along the path $P \rightarrow P' \rightarrow P''$. From $P$ to $P''$, there is no change in $\<x,\omega^\perp\>$, so there is no change in the term $e^{-\lambda\<x,\omega^\perp\>}$. So we only have to deal with the change of $f$.

Let us study the second figure. Denote $\varepsilon = p\tan\frac{\Delta\omega}{2}$. Then the length of $PP'$ is $(u+\varepsilon)\tan(\Delta\omega)$, this will multiply $\partial_\omega f(x)$, the directional derivative of $f(x)$ in the direction $\omega$. 

The length of $P'P''$ is $\frac{u+\varepsilon}{\cos\Delta\omega}+\varepsilon -u$, this will multiply $\partial_{\omega^\perp}f(x)$ approximately when $\Delta\omega$ is small.

To simplify these expressions, we can use their power series form:
$$\tan \Delta\omega = \Delta\omega + o(\Delta\omega) \;\;\mathrm{and}\;\; \frac{1}{\cos \Delta\omega} = 1 + o(\Delta\omega),$$ so
$$(u+\varepsilon)\tan(\Delta\omega) \approx (u+p\frac{\Delta\omega}{2})\,\Delta\omega = 
u\Delta\omega + o(\Delta\omega)$$ and
$$\frac{u+\varepsilon}{\cos\Delta\omega}+\varepsilon -u \approx 
\frac{u+p\frac{\Delta\omega}{2}}{1} + p\frac{\Delta\omega}{2} - u = p\Delta\omega.$$

So the change of $f$ between $P$ and $P''$ is altogether
$$u\Delta\omega\,\partial_\omega f(x) + p\Delta\omega\,\partial_{\omega^\perp}f(x),$$
where $u=\<x,\omega^\perp\>$.

The variable of integration in the original integral is $u$, the variable of integration after rotating the line by $\Delta\omega$ is $u'=u/\cos\Delta\omega=u(1+o(\Delta\omega))\approx u$. 

We divide the change of the integral by $\Delta\omega$ when differentiating with respect to $\omega$.
So the derivative w.r.t. the angle $\omega$ is
$$\partial_\omega R_\lambda f (\omega, p) = \int_{\<\omega,x\>=p} e^{-\lambda\<x,\omega^\perp\>} \,\left[ \<x,\omega^\perp\>\partial_\omega f(x) \,+\, p\,\partial_{\omega^\perp} f(x) \right] \, dx;$$
and it can be calculated, knowing the exponential Radon transform along lines avoiding $K$.

Using integration by parts for the second term of this expression,
\begin{equation}\label{eq:intbyparts}
\begin{split}
\int_{\<\omega,x\>=p} e^{-\lambda\<x,\omega^\perp\>} \, p\,\partial_{\omega^\perp} f(x) \, dx = \int_{-\infty}^\infty e^{-\lambda u} \, p\,\partial_u \{f(p\,\omega + u\,\omega^\perp)\} \, du =\\
= \left[e^{-\lambda u} \, p\,f(p\,\omega + u\,\omega^\perp)\right]_{-\infty}^\infty - \int_{-\infty}^\infty e^{-\lambda u}(-\lambda) \, p\,f(p\,\omega + u\,\omega^\perp) \, du =\\
= p\lambda \int_{\<\omega,x\>=p} e^{-\lambda\<x,\omega^\perp\>} \, f(x) \, dx = p\lambda R_\lambda f(\omega,p).
\end{split}
\end{equation}
This is a known quantity, so we will be able to subtract it from the derivative w.r.t. $\omega$. (Here $\left[e^{-\lambda u} \, p\,f(p\,\omega + u\omega^\perp)\right]_{-\infty}^\infty$ is zero because the exponential Radon transform is finite on $l(\omega,p)$.)
So
$$\partial_\omega R_\lambda f (\omega,p) = 
\int_{\<\omega,x\>=p} e^{-\lambda\<x,\omega^\perp\>}\<x,\omega^\perp\>\partial_\omega f(x) \, dx\; +\; p\lambda R_\lambda f(\omega,p).$$
This means, we can calculate $$\int_{\<\omega,x\>=p} e^{-\lambda\<x,\omega^\perp\>}\<x,\omega^\perp\>\partial_\omega f(x) \,dx.$$

Then using a similar argument as the one for differentiating $R_\lambda f(\omega,p)$ w.r.t. $p$ (partial integration of $\partial_{\omega} f(x)$ w.r.t. $p=\<\omega,x\>$ yields $f(x)$, and $\<x,\omega^\perp\>$ does not depend on $p$), we can see that integrating the above expression w.r.t. $p$ (renamed to $q$, from $p$ to $q_1$), and substituting $\omega_0$ for $\omega$, yields 
$$\int_p^{q_1} \int_{\<\omega_0,x\>=q} e^{-\lambda\<x,\omega_0^\perp\>}\<x,\omega_0^\perp\>\partial_\omega f(x) \,dx \; dq =$$
$$\int_{\<\omega_0,x\>=q_1} e^{-\lambda\<x,\omega_0^\perp\>}\<x,\omega_0^\perp\> f(x) \,dx - \int_{\<\omega_0,x\>=p} e^{-\lambda\<x,\omega_0^\perp\>}\<x,\omega_0^\perp\> f(x) \,dx,$$ where the first integral is given according to part e) of the system of conditions in $\mathcal{E}(\omega_0, \Omega, p_0)$, so the second integral can also be calculated for any $p>p_0$.

Now we are going to prove by induction that \emph{the values} $$R_\lambda^{(k)}f(\omega_0,p) := 
\int_{\<\omega_0,x\>=p} e^{-\lambda\<x,\omega_0^\perp\>}(\<x,\omega_0^\perp\>)^k \, f(x) \,dx$$
\emph{are all determined by the values of the exponential Radon transform 
 and are all finite, with the given system of conditions.}

Suppose we already know that the above statement holds for any $k$ such that $0\le k\le K$, for some $K$
($k\in\mathbb{N}$ and $K\in\mathbb{N}$). For $k=0$ it is obvious, and we have just proved it for $k=1$, so we have it for $K=1$. Now let's prove for a fixed $K\ge 1$ that consequently the statement is true for $k=K+1$, so for $0\le k\le K+1$, too.

Having the statement for $K$, differentiate $R_\lambda^{(K)}f(\omega,p)$ w.r.t. the angle $\omega$. As before, we move along a path that keeps 
$\<x,\omega^\perp\>$ fixed ($P\rightarrow P''$ in the second image).
So the differentiation yields
\begin{equation}
\begin{split}
\left.\partial_\omega R_\lambda^{(K)}f(\omega,p) \right|_{\omega=\omega_0}=\\
= \left.\int_{\<x,\omega\>=p} e^{-\lambda\<x,\omega^\perp\>}(\<x,\omega^\perp\>)^K \cdot \left[\<x,\omega^\perp\>\partial_\omega f(x) +  \, p \, \partial_{\omega^\perp}f(x) \right] \,dx\right|_{\omega=\omega_0}= \\
=\int_{\<x,\omega_0\>=p} e^{-\lambda\<x,\omega_0^\perp\>}(\<x,\omega_0^\perp\>)^{K+1}\partial_\omega f(x) \,dx \;+\\
+ \; p \int_{\<x,\omega_0\>=p} \left[ e^{-\lambda\<x,\omega_0^\perp\>}(-\lambda)(\<x,\omega_0^\perp\>)^K + e^{-\lambda\<x,\omega_0^\perp\>} K(\<x,\omega_0^\perp\>)^{K-1} \right] \, f(x) \,dx.
\end{split}
\end{equation}
Here we used integration by parts like in equation (~\ref{eq:intbyparts}): using the substitution $u=\<x,\omega_0^\perp\>$ and that $\partial_{\omega^\perp}=\partial_u$, and we also needed that 
$R_\lambda^{K}f(\omega_0,p)$ is finite that was included is the induction hypothesis.

In this last form, the first part of the second term is equal to $-\lambda p R_\lambda^{(K)}f(\omega_0,p)$ and the second part of the second term is equal to $Kp R_\lambda^{(K-1)}f(\omega_0,p)$. This means that the first term can be calculated as well. 

Integrating the first term w.r.t. $p$ as we did for $k=1$ (renamed to $q$, from $p$ to $q_{K+1}$) yields 
$$\int_p^{q_{K+1}} \int_{\<\omega_0,x\>=q} e^{-\lambda\<x,\omega_0^\perp\>}(\<x,\omega_0^\perp\>)^{K+1}\partial_\omega f(x) \,dx \; dq =$$
$$\int_{\<\omega_0,x\>=q_{K+1}}e^{-\lambda\<x,\omega_0^\perp\>}(\<x,\omega_0^\perp\>)^{K+1} f(x) dx - \int_{\<\omega_0,x\>=p} e^{-\lambda\<x,\omega_0^\perp\>}(\<x,\omega_0^\perp\>)^{K+1} f(x) dx$$ where the first integral is given according to part e) of the system of conditions in $\mathcal{E}(\omega_0, \Omega, p_0)$, so the second integral can also be calculated for any $p>p_0$. So we get that $R_\lambda^{K+1}f(\omega_0,p)$ can be calculated as well, along lines $l(\omega_0,p)$ avoiding the set $K$. The induction is completed.

\smallskip

Now we remind the reader for the following:
\begin{statement}\label{th:expdec} If a function $f:{\mathbb R}\to {\mathbb C}$ decayes exponentially 
and the integral of its product with every polynomial is 0, then the function is 0. 
\end{statement}

This can be proved using the uniqueness of the two-sided Laplace transform as follows. Although it is elementary and well known, we include the proof, because no good reference was found.

There are counter-examples to this statement if the decay is only faster than any polynomial~\cite{Helgason}. Our counter-examples for hyperbolic domains together with the proof of Theorem~\ref{th:supp2} also yield counter-examples for this statement with decay $e^{|x|^{1-\epsilon}}$ for any 
$0< \epsilon < 1$. We do not know of any other explicite counter-example. 

\emph{Proof:}
We take the two-sided Laplace transform ${\mathcal L}f(s) = \int^\infty_{-\infty} f(x) e^{-sx}\, dx$ of $f$. Because $f$ is of exponential decay, that is $f(x) e^{-\mu x}$ is bounded, for a fixed $\mu > 0$ the two-sided Laplace-transform  ${\mathcal L}f(s)$ is defined for any $0< s < \mu$ and complex analytic in the strip 
$0 < Re z < \mu$.  

Take the powes series expansion  of $e^{-sx}$ and interchange the summation and the integral: 
${\mathcal L}f(s) =  
\int^\infty_{-\infty} f(x)\Big|\sum_{n=0}^\infty \frac{(-1)^n}{n!}s^n x^n\Big]\, dx =
\sum_{n=0}^\infty \frac{(-1)^n}{n!}s^n \int^\infty_{-\infty}  f(x)x^n\, dx$. Our condition is that 
$\int^\infty_{-\infty}  f(x)x^n\, dx = 0$ for every $n \geq 0$, so the sum is 0, that is  ${\mathcal L}f(z) = 0$ for on the strip $0 < Re z < \mu$. According for example to~\cite{Widder} Theorem 6b, 
this implies that $f \equiv 0$.

To complete the proof, we only have to show that interchanging the sum and the integral is valid. 
We estimate the following two terms independently: 
${\mathcal L}f(s) = - \int^\infty_{0} f(-x) e^{sx}\, dx + \int^\infty_{0} f(x) e^{-sx}\, dx$. 

First, we estimate the tail of the first term. 

\begin{equation}
\begin{split}
 \Big|\int^\infty_{0} f(-x)\Big( \sum_{n=0}^N \frac{(sx)^n}{n!}- e^{sx} \Big)\, dx \Big| \leq 
C\int^\infty_{0} e^{-\mu x}\Big(e^{sx} - \sum_{n=0}^N \frac{(sx)^n}{n!}\Big)\, dx  = \\
C\int^\infty_{0} e^{(s-\mu)x} - C\sum_{n=0}^N \int^\infty_{0} e^{-\mu x} \frac{(sx)^n}{n!} =
\frac{C}{\mu - s} - \frac{C}{\mu} \sum_{n=0}^N \frac{s^n}{\mu^n} = \\
 \frac{C}{\mu}\Big( \frac{1}{1- \frac{s}{\mu}}  - \sum_{n=0}^N \frac{s^n}{\mu^n}  \Big).
\end{split}
\end{equation}
This is just the tail of the power series expansion of $\frac{C}{\mu}\frac{1}{1- \frac{s}{\mu}} $. 
It tends to 0. So the integral and the summation are interchangable.

Now we estimate the tail of the second term. This is an alternating series. So it can be estimated by 
its $(N+1)$th term.

\begin{equation}
\begin{split}
 \Big|\int^\infty_{0} f(x)\Big( \sum_{n=0}^N \frac{(sx)^n}{n!}- e^{-sx} \Big)\, dx \Big| \leq
C\int^\infty_{0} e^{-\mu x} \Big( \sum_{n=0}^N \frac{(sx)^n}{n!}- e^{-sx} \Big)\, dx \leq \\
C\int^\infty_{0}\int^\infty_{0}  e^{-\mu x} \frac{(sx)^{n+1}}{(n+1)!}\, dx = 
\frac{1}{\mu} \frac{s^{n+1}}{\mu^{n+1}}.
\end{split}
\end{equation}
So it tends to 0. So the integral and the summation is interchangable for this term too. $\square$

Observe that the statement above extends to any dimension. We can prove it by the classical Radon-transform on hyperplanes, because $\int^{\infty}_{-\infty} Rf(\omega, p) p^k\,dp = \int_{{\mathbb R}^n} f(x) \<x, \omega\>^k\,dx$. Of course, we can prove it directly too, but that is more technical.

\bigskip

On this point we need the conditions a) - d) from the definition of the class $\mathcal{E}(\omega_0, \Omega, p_0)$. They ensure that if $R_\lambda^{k}f(\omega,p)$ is zero for all $k\in\mathbb{N}$ and all $\omega\in\Omega$, $p>p_0$, then $f\equiv 0$ on $\mathcal{L}(\Omega,p_0)$. This implies that if there are two functions $f_1$ and $f_2$ with the same value of the exponential Radon transform on lines avoiding the set $K$, then because of the linearity of the exponential Radon transform, $f_1-f_2\equiv 0$, so they are identical.
\bigskip
$\square$

\medskip
\emph{Proof of Theorem~\ref{th:supp}.} 

All the requirements of Theorem~\ref{th:supp2}, apart from $a)$ and $c)$, are consequences of the exponential decay supposed in Theorem~\ref{th:supp}.
$a)$ and $c)$ are not implied because this simple decay condition does not ensure anything for continuity and the derivatives. But a standard technique: convolving $f$ belonging to this function class with a smooth function $g$ of compact support will help us show that Theorem~\ref{th:supp} still follows from Theorem~\ref{th:supp2}.
The condition e) will still be valid after the convolution, with $q_k=+\infty$.

Suppose we have a function $f$ that is not everywhere zero in $L^1_{loc}(K^c)$ and has exponential decay, but whose exponential Radon transform is zero along all straight lines in $K^c$ (that would be a counterexample for Theorem~\ref{th:supp}). We will take its convolution with a $\mathcal{C}^\infty$ function $g$ such that ${\rm supp}\, g$ is compact and ${\rm supp}\, g \subset B_{\frac{\varepsilon}{2}}(0)$ and the integral of $g$ is 1.
The convolution $f*g$ has exponential decay (inherited from $f$) and belongs to $\mathcal{C}^\infty$ (inherited from $g$).

Condition $e)$ is fulfilled by $f*g$ with $q_k=+\infty$, because $$\lim_{j\rightarrow\infty}\left(\int_{-\infty}^\infty u^k e^{-\lambda u}\, f(q_{kj}\, \omega_0 + u\, \omega_0^\perp)\, du\right)=0,$$ as a consequence of the exponential decay of $f$, and it is not spoiled by the convolution with $g$.

It is well known and we will prove a bit later that $$R_\lambda(f*g)(\omega,p)=(R_\lambda g*R_\lambda f)(\omega,p),$$ a one-dimensional convolution in the second variable. So $R_\lambda f(\omega,p)=0$ over all straight lines in $K^c$ implies  
$R_\lambda(f*g)(\omega,p)=0$ over all straight lines in $K_{\varepsilon}^c$, where $\varepsilon$ is bigger than the diameter of the support of $g$. 

But using the condition that $f$ is not identically zero in $L^1_{loc}(K^c)$, we can choose $g$ so that $f*g$ is not identically zero in $L^1_{loc}(K_{\varepsilon}^c)$, where $\varepsilon$ is bigger than the diameter of the support of $g$.
(The convolution $f*g$ of a function $f$ in $L^1_{loc}$ that is non-zero on a set of positive measure, with a non-zero non-negative continuous function $g$ of compact support, is not identically zero if the support of $g$  
is small enough.)

This leads to a contradiction with Theorem~\ref{th:supp2}, because $f*g$, restricted to $K_\varepsilon^c$, fulfils all the conditions of Theorem~\ref{th:supp2}, including $c)$ about the decay of the partial derivatives.
(We can think about $\varepsilon$ tending to zero, but it is not even necessary, we just have to choose it in such a way that $f*g$ is not zero.)

This way we proved that $R_\lambda f(\omega,p)=0$ over all straight lines in $K^c$ implies  $f = 0$ in $L^1_{loc}(K^c)$. 

\bigskip
For the proof of $R_\lambda(f*g)(\omega,p)=(R_\lambda g*R_\lambda f)(\omega,p)$,
let us denote $f*g$, the (two-dimensional) convolution of $f$ and $g$, by $h$: $h(x)=\int_{\mathbb{R}^2}\,f(x-y)g(y)\,dy$. Then
$$R_\lambda h(\omega,p) = \int_{\<\omega,x\>=p}\,e^{-\lambda\<\omega^\perp,x\>}\,h(x)\,dx=$$ $$=\int_{\<\omega,x\>=p}\,e^{-\lambda\<\omega^\perp,x\>}\,\int_{\mathbb{R}^2}\,f(x-y)g(y)\,dy\,dx=$$
$$=\int_{\<\omega,x\>=p}\,e^{-\lambda\<\omega^\perp,x-y\>}e^{-\lambda\<\omega^\perp,y\>}\,\int_{\mathbb{R}^2}\,f(x-y)g(y)\,dy\,dx=$$
$$=\int_{\mathbb{R}^2}\,e^{-\lambda\<\omega^\perp,y\>}\,g(y)\left\{\int_{\<\omega,x-y\>=p-\<\omega,y\>}\,f(x-y)\,e^{-\lambda\<\omega^\perp,x-y\>}\,d(x-y)\right\}dy$$
where $d(x-y)=dx$ for a fixed $y$. Let us introduce the notation $\<\omega,y\>=q$, then we can write the above expression as
$$=\int_{\mathbb{R}^2}\,e^{-\lambda\<\omega^\perp,y\>}\,g(y)\,R_\lambda f(\omega,p-q)\,dy,$$
where $q$ is dependent on $y$. Here $\omega$ and $p$ are fixed, so for a fixed $q$, $R_\lambda f(\omega,p-q)$ is a constant. We integrate in $\mathbb{R}^2$, that can be decomposed to lines perpendicular to the vector $\omega$; and we can integrate along these lines first, then with respect to the signed distance  of the line from the origin (which is $q$), obtaining that the above expression equals
$$=\int_{\mathbb{R}}\,R_\lambda g(\omega,q)\,R_\lambda f(\omega,p-q)\,dq=(R_\lambda g*R_\lambda f)(\omega,p),$$
a one-dimensional convolution.

$\square$

\section{Further statements}

Using the Helgason moment-condition~\cite{Helgason}, that is $\int_{\mathbb R} Rf(\omega, p) p^k\,dp$ is an at most $k$th degree polynomial of $\omega\in S^{n-1}$, provided $f$ decays faster than any polynomial, 
we can easily derive the following well known statement. 

Let $\Omega\subset{\mathbb S^1}$ an infinite set, and $f: \;{\mathbb R^2}\rightarrow {\mathbb R}$ a continuous function with compact support. Assume that $Rf(\omega, p)=0$ for every $\omega\in\Omega$ and $p\in\mathbb{R}$. Then $f\equiv 0$. Note that the condition that $f$ is decaying faster than any 
polynomial is not enough, even if $\Omega = S^{n-1}$. But by Statement~\ref{th:expdec} we can apply the proof for functions with exponential decay too. This version is not mentioned in the literature to our knowledge.

In ${\mathbb R^n}$ the same is true, but $\Omega\subset {\mathbb S^{n-1}}$ has to be Zariski dense! 
(If a homogenuous polynomial $\equiv 0$ on $\Omega$ then it is zero.)

We can extend this well known
result as follows. The difference is, that we have a kind of obstacles. This contains the case of 
limited-angle tomography with obstacles. This might be of interest from practical point of view.  

For the presence of the obstacle we have to pay by slightly stronger conditions on $\Omega \subseteq S^{n-1}$.
\begin{theorem}\label{th:supp3}
Let $f: \;\mathbb R^2\rightarrow\mathbb R$ be a continuous function and $K$ is a convex, closed set.

Let $\Omega\subset\mathbb S^1$ be a perfect set (every point is a limit point). 
Let $K^\#$ be the union of those lines which avoid $K$ and their normal vectors are in $\Omega$.

Assume that $f$ admits the decay conditions of Theorem~\ref{th:supp2} for every $\omega_0 \in \Omega$ and 
$R_\lambda f(\omega,p)=0$ if $l(\omega,p)\bigcap K=\emptyset$ and $\omega\in\Omega$. 
Then $f\equiv 0$ on $K^\#$.
\end{theorem}

\emph{Proof:} The immediate consequence of the proof of Theorem~\ref{th:supp2}. $\square$.

To our knowledge this statement is new even for $\lambda=0$, the case of classical Radon transform.

Note that for the exponential Radon transform there is no moment condition, so even the $K=\emptyset$ case is non-trivial. Also note that the Helgason moment condition~\cite{Helgason} cannot be applied if 
$K\neq \emptyset$.

We can extend this to higher dimensions for the classical Radon transform if $\Omega$ is rich enough.

By proof mining of Theorem~\ref{th:supp2}, we can easily find sufficient conditions for $\Omega$ even in the $n=2$ case, that is  Theorem~\ref{th:supp3}.

Such a sufficient condition is as follows.

Let $\omega\in\Omega$. We say that $\iota\in\mathbb S_x^{n-2}$ (the $\mathbb S_\omega^{n-1}$ sphere in the tangent space $T_\omega\mathbb S^{n-1}$) is a limit direction if there is a sequence of points $\omega_1, ..., \omega_k,...\rightarrow\omega$ such that $\frac{\omega_k-\omega}{|\omega_k-\omega|}\rightarrow\iota$.
The set of all limit directions for $\omega\in\Omega$ be $I(\omega)$. 

\begin{theorem}\label{th:supp4}
If $\Omega\subset\mathbb S^{n-1}$ is perfect and $I(\omega)$ is Zariski dense for every $\omega\in\Omega$ then if $K$ is closed  and $f:\;\mathbb R^n\rightarrow\mathbb R$ satisfies the multi-dimensional version of the decay conditions of Theorem~\ref{th:supp2} for every $\omega_0\in \Omega$, and $Rf(\Omega,p)=0$ for every hyperplane $H(\omega,p)$ avoiding $K$ and $\omega\in\Omega$ then $f\equiv 0$ on $K^\#$.
\end{theorem}

\emph{Proof:} Apply the same technique as in Theorem~\ref{th:supp2} for hyperplanes almost perpendicular (in  second order) to an apropriate two-plane 
$\theta(\omega_0, \eta)$ spanned by the vectors $\omega_0 $ and $\eta\perp \omega_0$. Change this plane at every step of the derivations. So we get the 
integrals $\int_{\<x,\omega_0\> = p} f(x) \<x, \eta_1\>\cdots \<x, \eta_k\>\,dx$. This implies that the integral $\int_{\<x,\omega_0\> = p} f(x)p(x) =0$
for any polynomial $p(x)$ defined on the points of the hyperplane $H(\omega_0, p) = \{x\in {\mathbb R}^n\}: \<x,\omega_0\> = p \}$. Then apply Statement~\ref{th:expdec}. $\square$.

Note that both in Theorem~\ref{th:supp3} and Theorem~\ref{th:supp4} functions with exponential decay automatically satisfy the decay conditions.
 
The generalization of Theorem~\ref{th:supp2} is as follows. Let $R$ be the classical Radon-transform on ${\mathbb R}^n$, that is integration on hyperplanes, and let $R_{n-2}$ be the integration on $n-2$ planes. Here $\Omega\subseteq S^{n-1}$ is an open set containing $\omega_0$.
\begin{theorem}\label{th:supp5}
If  $f\in \mathcal E(\omega_0, \Omega,  p_0)$ 
and $Rf(\xi) = 0$ for every hyper-plane $\xi(\omega, p)$, for $\omega \in \Omega$ and $p\geq p_0$, then $f = 0$.
\end{theorem}
\emph{Proof:} 
Apply the same technique as in Theorem~\ref{th:supp2} for hyperplanes perpendicular to $\omega_0$ and a fixed $\eta\in S^{n-1}$ and $\eta\perp\omega_0$. By applying Fubini's theorem and Theorem~\ref{th:supp2} we get that $R_{n-2}f(\zeta_{n-2}) = 0$ for any $(n-2)$-plane $\zeta_{n-2}$ which is perpendicular to $\omega_0$ and in a hyperplane $\xi(\omega_0, p)$ whose oriented distance $p$ from $O$ is at least $p_0$. This means that the Radon transform of $f$ is zero in every such hyperplane $\xi(\omega_0, p)$. 
Note that in this case the Radon transform integrates on $(n-2)$-planes, which are hyperplanes in the hyperplane $\xi(\omega, p)$. Also note that the 
$(n-2)$-plane transform $R_{n-2}$ does not affect the decay conditions. $\square$


To make things simple, and easier to remember, we state it by stronger decay conditions too, as in Theorem~\ref{th:supp}.
\begin{theorem}\label{th:supp6}
If $K \subset {\mathbb R}^n$ is a convex closed set, and $f:K^c\to \mathbb R$ is a continuous function with exponential decay 
and $Rf(\xi) = 0$ for every hyper-plane $\xi$ avoiding $K$, that is in $\xi \subset K^c$, then $f = 0$.
\end{theorem}

Note that Theorem~\ref{th:supp2} and Theorem~\ref{th:supp5} is not a consequence of Theorem~\ref{th:supp3} and Theorem~\ref{th:supp4}, respectively. This is  because the decay conditions of Theorem~\ref{th:supp2} are required only for one direction $\omega_0$, not for all the directions $\omega_0\in \Omega$ as in the latter two theorems.   

The Theorem~\ref{th:supp6} above is also close to be sharp, as in the plane. We can prove it by the rotation-symmetric extension of our planar counter-examples to ${\mathbb R}^n$. But the geometry is a bit more complicated. We say that a closed convex set $K$ is strongly hyperbolic, if it does not contain a whole line, but contains $n$ independent directions of half-lines. We say that $K$ is parabolic, if it does not contain a whole line,  there are no $n$ independent half-lines, but there is at least one. For strongly hyperbolic  domains we have the rotation symmetric counter-example which decays as 
$e^{-|x|^{1-\epsilon}}$, while for parabolic ones decays as $e^{-|x|^{1/2-\epsilon}}$ for a fixed $\epsilon > 0$.
The easy details are left to the reader.

Let $\omega_0$ be as before, $\lambda > 0$ and $a_0(x) = \exp(-\lambda\cdot\<x, \omega^\perp_0\>)$. Then for $\<\omega^\perp, \omega^\perp_0\> > 0$ the attenuated Radon transform $R_a$, as defined in~\cite{Novikov},~\cite{Natterer}, are the same for all such $\omega$'s up to a positive multiplicative function.  

\begin{question}
Can we extend our results to attenuated Radon transforms with sufficiently smooth weights '$a$' which are exponentially close to '$a_0$'? That is $|a(x) - a_0(x)| <C\cdot\exp(-\mu\cdot|x|)$, where $C> 0$ is a suitable constant and $\mu > \lambda > 0$.
\end{question}

\end{document}